\documentclass[12pt]{article}
\usepackage{amssymb}

\input{tcilatex}

\begin{document}

\author{Naichung Conan Leung\thanks{%
This research is supported by NSF/DMS-9626689}}
\title{Uniformization of four manifolds }
\date{May 6, 1997}
\maketitle

\begin{abstract}
We characterize those Einstein four manifolds which are locally symmetric
spaces of noncompact type. Namely they are four manifolds which admit
solutions to the (non-Abelian) Seiberg Witten equations and satisfy certain
characteristic number equality.
\end{abstract}

\section{Introduction}

The main theorem in this paper is:

\begin{theorem}
Let $X$ be an Einstein four manifold with $\sigma \left( X\right) =0.$
Suppose for some $Spin^c$ structure $L$ on $X$, the non-Abelian Seiberg
Witten equations 
\begin{eqnarray*}
F_A^{+} &=&\tau \left( \phi \right)  \\
\frak{D}_A\phi  &=&0,
\end{eqnarray*}
on $S_L^{+}$ has a solution, then the universal covering of $X$ must be $%
\Bbb{R}^4,B^2\times B^2$ or $B^4$ with a locally symmetric Riemannian metric.
\end{theorem}

Combining with the work of Yau \cite{Y} in the K\"{a}hler surface case and
LeBrun \cite{LB} in the general case for the complex ball quotient $\Gamma
\backslash B_{\Bbb{C}}^2$, we can characterize every four dimensional
Einstein manifold $X$ which is locally symmetric of noncompact type by using
(non-Abelian) Seiberg Witten equations. Recall that in four dimension, the
following is a complete list of symmetric spaces of noncompact type: 
\[
\begin{array}{lll}
B^4, & B_{\Bbb{C}}^2, & B^2\times B^2.
\end{array}
\]
Therefore we obtain:

\begin{theorem}
Let $X$ be an Einstein four manifold with $\chi \left( X\right) >0$. Then $X$
is a locally symmetric space of noncompact type iff one of the followings
holds

(1) the non-Abelian Seiberg Witten equation for $S_L^{+}$ has a solution and 
$\sigma \left( X\right) =0$.

(2) the Seiberg Witten equation for $L$ has a solution and $3\sigma \left(
X\right) =\chi \left( X\right) .$

The universal cover of $X$ is isometric to $B^2\times B^2$ or $B^4$ in (1)
and $B_{\Bbb{C}}^2$ in (2).
\end{theorem}

\begin{remark}
Kotschick points out to the author that any Einstein metric on a four
manifold with $\chi \left( X\right) =0$ must be flat.
\end{remark}

\begin{remark}
It is desirable to define Donaldson type invariant using non-Abelian Seiberg
Witten equations. Because one would be able to replace the existence of
Seiberg Witten solution by a (differential) topological condition. In
particular, it will imply that Einstein metric on four dimension locally
symmetric space of noncompact type is \textit{unique}. To the author's
knowledge, uniqueness of Einstein metric on $S^4$ is still open.
\end{remark}

In dimension two and three, Einstein metric has constant sectional
curvature. Namely their universal covers are isometric to $S^{n},\Bbb{R}^{n}$
or $B^{n}$ where $B^{n}=\left\{ x\in \Bbb{R}^{n}:\left| x\right| <1\right\} $
is the unit ball with the hyperbolic metric. In particular they coincide
with locally symmetric spaces in these dimensions.

In dimension four Einstein metric may not be locally symmetric, for instance
the K3 surface with Yau's metric \cite{Y} is not locally symmetric. Hitchin 
\cite{H} showed that there is topological obstruction to the existence of
Einstein metric in this dimension. We combine Hitchin's original ideas and
the use of (non-Abelian) Seiberg Witten equations to determine which
Einstein metric is locally symmetric of noncompact type. Recall that in
dimension four, there are only three symmetric spaces $G/K$ of noncompact
type, namely

(i) the real four ball 
\[
SO^0\left( 4,1\right) /SO\left( 4\right) =B^4, 
\]
with the real hyperbolic metric;

(ii) the complex two ball 
\[
SU\left( 2,1\right) /U\left( 2\right) =B_{\Bbb{C}}^2=\left\{ x\in \Bbb{C}%
^k:\left| x\right| <1\right\} , 
\]
with the complex hyperbolic metric and

(iii) the reducible case 
\[
SO^0\left( 2,1\right) \times SO^0\left( 2,1\right) /SO\left( 2\right) \times
SO\left( 2\right) =B^2\times B^2, 
\]
which is the product of the unit disks.

Each of these spaces carries canonical (non-unitary) flat connection which
provides a solution to the (non-Abelian) Seiberg Witten equations. We will
reverse the process and use a solution of the (non-Abelian) Seiberg Witten
equations to construct a perturbed non-unitary anti-self-dual connection
over $X.$ Existence of such connection implies certain characteristic number
inequality for $X.$ When the equality holds (which depends only on the
homotopy structure of $X$) then such connection is in fact projectively
flat. By analyzing the flatness condition, we will prove in section four
that these informations are sufficient to show that $X$ must be a locally
symmetric space of noncompact type.

The pioneer work in this direction is due to Hitchin \cite{H} who showed
that an Einstein four manifold satisfies the following characteristic number
inequality 
\[
\left| \sigma \left( X\right) \right| \leq \frac 23\chi \left( X\right) , 
\]
where $\sigma \left( X\right) =b_{+}^2\left( X\right) -b_{-}^2\left(
X\right) $ is the signature of $X$ and $\chi \left( X\right) $ is the Euler
characteristic of $X.$ These are the only characteristic numbers in
dimension four and they depend solely on the homotopic type of the space $X.$

Moreover Hitchin showed that $\left| \sigma \left( X\right) \right| =\frac
23\chi \left( X\right) $ if and only if $X$ is flat or its universal cover
is a K3 surface. They can be distinguished by their Euler characteristics
being zero or not.

When the Einstein metric is also K\"{a}hler then the above characteristic
number inequality can be strengthened to 
\[
\sigma \left( X\right) \leq \frac 13\chi \left( X\right) , 
\]
or equivalently 
\[
c_1^2\left( X\right) \leq 3c_2\left( X\right) . 
\]

In \cite{Y} Yau showed that any K\"{a}hler surface with negative Ricci
curvature carries a unique K\"{a}hler Einstein metric. Therefore we have $%
c_1^2\left( X\right) \leq 3c_2\left( X\right) $ for any such K\"{a}hler
surface and when the equality sign holds then the universal cover of $X$ is
the complex two ball $B_{\Bbb{C}}^2.$

After the introduction of Seiberg Witten theory \cite{W}, LeBrun \cite{LB}
showed that in order to get the characteristic number inequality $\sigma
\left( X\right) \leq \frac 13\chi \left( X\right) $ the K\"{a}hlerian
condition can be replaced by the existence of solution to the Seiberg Witten
equation. Therefore Einstein four manifold with non-trivial Seiberg Witten
invariant and $\sigma \left( X\right) =\frac 13\chi \left( X\right) $ is
covered by $B_{\Bbb{C}}^2.$ In \cite{L} the author interpreted LeBrun work
as a construction of certain perturbed $U\left( 2,1\right) $ anti-self-dual
connection on the bundle $S_L^{\Bbb{-}}\bigoplus L$ (where $L$ defines the $%
Spin^c$ structure on $X$) which will be projectively flat if $\sigma \left(
X\right) =\frac 13\chi \left( X\right) .$

Notice that the number of solutions to the Seiberg Witten equations counted
with multiplicity is a differentiable invariant \cite{W}. As a result of it,
we can show that Einstein metric on $X=\Gamma \backslash B_{\Bbb{C}}^2$ is
unique up to diffeomorphisms.

In \cite{L} the author showed that if $X$ has a K\"{a}hler Einstein metric
with $c_1\left( X\right) <0$ such that both orientations has non-trivial
Seiberg Witten invariant then 
\[
\sigma \left( X\right) \geq 0. 
\]
Moreover the equality sign holds if and only if $X$ is covered by $B^2\times
B^2.$ Later \cite{K} generalized this result to any complex surfaces of
general type.

We are going to further develop these ideas to study four dimensional
locally symmetric spaces in general which include the quotient of the real
ball $B^4=SO^0\left( 4,1\right) /SO\left( 4\right) $. Instead of trying to
reconstruct the $SO^0\left( 4,1\right) $ flat connection, we shall construct
a $Sp\left( 1,1\right) $ connection on some rank four complex vector bundle (%
$Sp\left( 1,1\right) $ is the universal (double) cover of $SO^0\left(
4,1\right) $). It turns out that the correct bundle to look at is $%
S^{+}\bigoplus S^{-}$ and we need to study the non-Abelian Seiberg Witten
equations on $E=S^{+}.$ If $X$ does not have a spin structure, then we shall
use a $Spin^c$ structure which always exists on a four manifold. This
approach can be used to deal with both $B^4$ and $B^2\times B^2$ cases
together. For $B^2\times B^2,$ the canonical flat connection will live in a
rank four complex vector bundle via the natural embedding $U\left(
1,1\right) \times U\left( 1,1\right) \subset U\left( 2,2\right) $. The rank
four bundle in question is again $S^{+}\bigoplus S^{-}.$ In this approach,
we do not need to impose any K\"{a}hlerian assumption on $X.$

\section{Locally symmetric space in dimension four}

As explained in last section that there are only three kinds of locally
symmetric space of noncompact type in dimension four, namely their universal
covering $\tilde{X}$ is $B^{4},B_{\Bbb{C}}^{2}$ or $B^{2}\times B^{2}.$

When $\tilde{X}=B^{4}$ or $B^{2}\times B^{2}$ we have $\sigma \left(
X\right) =0$ by Hirzebruch proportionality principle because their compact
duals are $S^{4}$ and $S^{2}\times S^{2}$ and they both have zero signature.
We also have $\chi \left( X\right) >0.$

When $\tilde{X}=B_{\Bbb{C}}^{2}$ we have $3\sigma \left( X\right) =\chi
\left( X\right) >0$ and there is a canonical representation $\pi _{1}\left(
X\right) \rightarrow \Bbb{P}U\left( 2,1\right) $ up to conjugation.

When $\tilde{X}=B^{2}\times B^{2}$ then we have a projective representation
of $\pi _{1}\left( X\right) $ to $\Bbb{P}U\left( 1,1\right) \times \Bbb{P}%
U\left( 1,1\right) .$ When $\tilde{X}=B^{4}$ we have a projective
representation of $\pi _{1}\left( X\right) $ to $SO^{0}\left( 4,1\right) .$
Instead of finding a flat connection on a real rank five vector bundle over $%
X,$ we will construct a complex rank four vector bundle with a projectively $%
U\left( 2,2\right) $ connection and show that the image of its holonomy does
lie inside $Sp\left( 1,1\right) \subset U\left( 2,2\right) .$

Observe that we have a canonical isomorphism of Lie groups: 
\[
Spin^{0}\left( 4,1\right) =Sp\left( 1,1\right) 
\]
extending the well-known isomorphism 
\[
Spin\left( 4\right) =Sp\left( 1\right) \times Sp\left( 1\right) . 
\]
They give an isomorphism between symmetric spaces: 
\[
SO^{0}\left( 4,1\right) /SO\left( 4\right) =Spin^{0}\left( 4,1\right)
/Spin\left( 4\right) =Sp\left( 1,1\right) /Sp\left( 1\right) \times Sp\left(
1\right) , 
\]
which is the real four ball $B^{4}.$

On $B^4$ there is a flat $SO^0\left( 4,1\right) $ bundle $T\bigoplus \Bbb{R}%
. $ The corresponding $Sp\left( 1,1\right) $ flat bundle is $S^{+}\bigoplus
S^{-}.$ In particular there is an anti-self-dual connection on $%
S^{+}\bigoplus S^{-}$ which restricts to the Levi-Civita connection on $%
S^{-}.$ For compact $X=\Gamma \backslash B^4$ with $\Gamma $ a torsion free
subgroup of $SO^0\left( 4,1\right) .$ We have $\pi _1\left( X\right) =\Gamma 
$ and the above mentioned $SO^0\left( 4,1\right) $ flat connection will
descend to give a $SO^0\left( 4,1\right) $ flat bundle over $X.$

Instead of the representation of $\pi _1\left( X\right) $ into $SO^0\left(
4,1\right) $ or $SO^0\left( 2,1\right) \times SO^0\left( 2,1\right) ,$ we
want to lift these representations to $Spin^0\left( 4,1\right) $ and $%
Spin^0\left( 2,1\right) \times Spin^0\left( 2,1\right) $ and use the
identifications $Spin^0\left( 4,1\right) =Sp\left( 1,1\right) $ and $%
Spin^0\left( 2,1\right) =SU\left( 1,1\right) .$ Both $Sp\left( 1,1\right) $
and $SU\left( 1,1\right) \times SU\left( 1,1\right) $ are closed subgroup of 
$U\left( 2,2\right) $ and therefore they give us a flat connection on a rank
four complex vector bundle over $X.$ The bundle in question is $%
S^{+}\bigoplus S^{-}.$

In general we may not be able to lift it to a $Spin^0$ flat connection. To
resolve this problem, we introduce a $Spin^c$ structure on $X$ and study 
\textit{projectively flat} $U\left( 2,2\right) $ connection over $X.$ Recall
that $Spin^c\left( 4\right) =Spin\left( 4\right) \times _{\Bbb{Z}_2}U\left(
1\right) $ and there is a short exact sequence 
\[
0\rightarrow \Bbb{Z}_2\rightarrow Spin^c\left( 4\right) \rightarrow SO\left(
4\right) \times U\left( 1\right) \rightarrow 0. 
\]
Therefore any $Spin^c$ structure on $X$ induces a $U\left( 1\right) $ line
bundle $L$ over $X.$ The homomorphism to $SO\left( 4\right) $ just gives us
back the tangent bundle of $X.$ Using the fact that $Spin^c\left( 4\right)
=SU\left( 2\right) \times SU\left( 2\right) ,$ we have two homomorphisms $%
Spin^c\left( 4\right) \rightarrow U\left( 2\right) .$ They induce two rank
two complex vector bundles over $X$ which are called the positive and
negative spinor bundles and we denote them by $S_L^{+}$ and $S_L^{-}.$
Formally they are $S^{+}\bigotimes L^{1/2}$ and $S^{-}\bigotimes L^{1/2}.$

Now we assume that the universal cover $\widetilde{X}$ of $X$ is either $B^4$
or $B^2\times B^2$. Then $\widetilde{X}$ has a canonical flat connection on $%
S^{+}\bigoplus S^{-}$ which may not descend to one on $X.$ After choosing a $%
Spin^c$ structure on $X,$ we consider the rank four complex vector bundle $%
S_L^{+}\bigoplus S_L^{-}$ over $X$ and we fix any Hermitian connection $D_L$
on $L,$ then the canonical flat connection on $S^{+}\bigoplus S^{-}$ twisted
by $D_L$ can be descended to $X$ and give a projectively flat connection
over $X.$ From our later discussions we shall see that such a projectively
flat $U\left( 2,2\right) $ connection provides a perturbed anti-self-dual
connection over $X$ up to gauge transformations. This perturbed
anti-self-dual connection will be constructed from a solution to the Seiberg
Witten equations and using the Einstein metric on $X.$

\section{(non-Abelian) Seiberg Witten equations}

In this section we assume that $X$ is a four manifold with a Riemannian
metric $g$ and a $Spin^c$ structure $L,$ that is $c_1\left( L\right) \equiv
w_2\left( X\right) \left( \mbox{mod}2\right) .$ We fix a connection $D_L$ on 
$L$ throughout our discussions. Before introducing the non-Abelian Seiberg
Witten equations, let us first review some well-known properties of
curvature decomposition in dimension four.

From the choice of a $Spin^c$ structure, we have the corresponding positive
and negative spinor bundle $S_L^{+}$ and $S_L^{-}$, they have canonical
connections induced from the Levi-Civita connection of $X$ and $D_L$ on $L.$
We now explain the relationship between the curvature of $\Lambda _{\pm
}^2,S_L^{\pm }$ and the Riemann curvature tensor of $X.$ The Riemann
curvature tensor of $X$ defines a curvature operator $Rm$ on the space of
two forms $\Lambda ^2$ on $X$ and the Hodge star operator $*$ of $X$
decomposes $\Lambda ^2$ into $\pm 1$ eigenspaces, namely the space of
self-dual and anti-self-dual forms: 
\[
\Lambda ^2=\Lambda _{+}^2\bigoplus \Lambda _{-}^2. 
\]
With respect to this decomposition $Rm$ has the following blocks
decomposition: 
\[
Rm=\left( 
\begin{array}{ll}
A & B \\ 
B^{*} & D
\end{array}
\right) , 
\]
where $A=A^{*}\in End\left( \Lambda _{+}^2\right) $, $D=D^{*}\in End\left(
\Lambda _{-}^2\right) $ and $B\in Hom\left( \Lambda _{-}^2,\Lambda
_{+}^2\right) .$ If we denote the scalar curvature of $X$ by $s$ then $%
s=4TrA=4TrD.$ Moreover $W^{+}=A-\frac 13TrA$ and $W^{-}=D-\frac 13TrD$ is
called the self-dual and anti-self-dual Weyl curvature of $X.$ The sum $%
W=W^{+}+W^{-}$ is called the Weyl curvature which depends only on the
conformal class of $g.$ In fact $W=0$ if and only if $X$ is a conformally
flat manifold.

We also have $B=Rc^0$ the trace free part of the Ricci curvature tensor.

By the Hirzebruch signature theorem, we have

\[
\sigma \left( X\right) =\frac 1{12\pi ^2}\int_M\left( \left| W^{+}\right|
^2-\left| W^{-}\right| ^2\right) dv_g. 
\]
In particular, if $X$ has zero signature, then $W^{-}=0$ if and only if $%
W=0. $ That is $X$ is a conformally flat manifold.

There are also natural isomorphisms 
\begin{eqnarray*}
\Lambda _{+}^2 &=&\mathbf{su}\left( S_L^{+}\right) \\
\Lambda _{-}^2 &=&\mathbf{su}\left( S_L^{-}\right) .
\end{eqnarray*}

Atiyah, Hitchin and Singer \cite{AHS} observed that $X$ is Einstein if and
only if the Levi-Civita connection on $\Lambda _{-}^2$ is anti-self-dual
connection. They also showed that $X$ is Einstein if and only if the trace
free part of the curvature on $S_L^{-}$ is anti-self-dual. Moreover $%
W^{-}=0=s$ if and only if the trace free part of the curvature on $S_L^{-}$
is self-dual. Corresponding statements for $\Lambda _{+}^2$ and $S_L^{+}$
also hold if we interchange self-dual by anti-self-dual.

In particular if $X$ is an Einstein manifold with $\sigma \left( X\right) =0$
such that the trace free part of the curvature on $S_L^{-}$ is self-dual
then the Riemann curvature tensor of $g$ vanishes identically. That is $X$
is a flat manifold.

Now we are going to introduce the (non-Abelian) Seiberg Witten equation: let 
$E$ be a complex Hermitian vector bundle over $X$ of rank $r$. Suppose $D_A$
is any Hermitian connection on $E$ and $\phi \in \Gamma \left( X,Hom\left(
S_L^{+},E\right) \right) =\Gamma \left( X,S_{L^{-1}}^{+}\bigotimes E\right)
. $ We consider $\phi \bigotimes \bar{\phi}\in \Gamma \left( X,\mathbf{u}%
\left( S_L^{+}\right) \bigotimes End\left( E\right) \right) .$ Under the
identification $\mathbf{su}\left( S_L^{+}\right) =\Lambda _{+}^2$, we
project $\phi \bigotimes \bar{\phi}$ to the trace free part of $\mathbf{u}%
\left( S_L^{+}\right) $ and denote its image by $\tau \left( \phi \right)
\in \Omega _{+}^2\left( X,End\left( E\right) \right) .$

We consider the following non-Abelian Seiberg Witten equations 
\begin{eqnarray*}
F_{A}^{+} &=&\tau \left( \phi \right) \\
\frak{D}_{A}\phi &=&0.
\end{eqnarray*}
Here $\frak{D}_{A}$ is the Dirac operator on $S_{L^{-1}}^{+}$ twisted by the
connection $D_{A}$ on $E.$

These equations are natural generalization of Seiberg Witten equations \cite
{W} to higher rank bundles. The Seiberg Witten equations have had numerus
applications in four manifold geometry, for examples \cite{FS}, \cite{FM}, 
\cite{KM}, \cite{LL}, \cite{MST}, \cite{T} and many others. Their higher
rank analog has also been studied extensively. For example, Pidstrigach,
Tyurin \cite{P}, \cite{PT} and also Feehan, Leness \cite{FL} use these
equations in an attempt to prove the equivalency of the Donaldson invariant
and the Seiberg Witten invariant with much progress. In the context of
Kahler surfaces, they are studied by Bradlow, Garcia-Prada \cite{BG} and
Okonek, Teleman \cite{OT}.

In this paper we generalize the arguments in \cite{L} to relate solutions to
these equations and perturbed anti-self-dual connection. This is crucial for
producing non-unitary flat connections over $X.$ Most of our discussions for
the rest of this section are based on methods used in section four of \cite
{L}. Even though almost everything is the same for general rank $r$ vector
bundle, we shall restrict our attention to the case when $r=2$ and $\det E=L$
for our purposes.

Now we suppose that $\left( D_A,\phi \right) $ is a solution to the
non-Abelian Seiberg Witten equations. Via the Clifford multiplication, there
is a canonical homomorphism 
\begin{eqnarray*}
\Gamma \left( X,Hom\left( S_L^{+},E\right) \right) &\rightarrow &\Omega
^1\left( X,Hom\left( S_L^{-},E\right) \right) \\
\phi &\rightarrow &\tilde{\phi}.
\end{eqnarray*}
From \cite{L}, we have 
\[
P_{-}\left( D_A\tilde{\phi}\right) =0, 
\]
for any harmonic spinor $\phi .$ Here $P_{-}$ is the orthogonal projection
to the anti-self-dual two forms. Moreover by projecting to the self-dual two
forms part we have 
\begin{eqnarray*}
P_{+}\left( \tilde{\phi}\tilde{\phi}^{*}\right) &=&Tr_E\tau \left( \phi
\right) I_{S_L^{-}}\in \Omega _{+}^2\left( \mathbf{u}\left( S_L^{-}\right)
\right) , \\
P_{+}\left( \tilde{\phi}^{*}\tilde{\phi}\right) &=&-2\tau \left( \phi
\right) \in \Omega _{+}^2\left( \mathbf{u}\left( E\right) \right) .
\end{eqnarray*}

We consider the connection $\Bbb{D}_{\left( A,\phi \right) }=\left( 
\begin{array}{ll}
\nabla _{L}^{LC} & \tilde{\phi} \\ 
\tilde{\phi}^{*} & D_{A}
\end{array}
\right) $ on $\Bbb{E}=S_{L}^{-}\bigoplus E$. $\Bbb{D}_{\left( A,\phi \right)
}$ preserves the Hermitian form $h_{S_{L}^{-}}-h_{E}$ on $\Bbb{E}$ and
therefore is a $U\left( 2,2\right) $ connection. If we extend the Hodge star
operator $*$ on $X$ to $\Omega ^{2}\left( X,End\left( \Bbb{E}\right) \right) 
$ by 
\[
\ast _{E}\left( 
\begin{array}{ll}
A & B \\ 
C & D
\end{array}
\right) =\left( 
\begin{array}{ll}
\ast A & -*B \\ 
-*C & *D
\end{array}
\right) . 
\]
Using $*_{E}$ we can decompose $\Omega ^{2}\left( X,End\left( \Bbb{E}\right)
\right) $ into self-dual and anti-self-dual parts.

Now if we decompose the curvature $\Bbb{F}_{\left( A,\phi \right) }$ of $%
\Bbb{D}_{\left( A,\phi \right) }$ into self-dual and anti-self-dual
accordingly, we obtain 
\[
\Bbb{F}_{\left( A,\phi \right) }^{+}=\left( 
\begin{array}{ll}
Rc^0-F_L^{+} & 0 \\ 
0 & -F_L^{+}
\end{array}
\right) . 
\]
Here $F_L$ is the curvature of $D_L.$ The proof of this formula is by direct
computations and it is identical to the rank one case treated in section
four of \cite{L}.

Hence we have shown the following result:

\begin{proposition}
Let $X$ be an Einstein four manifold with $Spin^c$ structure $L$ on $X$.
Suppose that the non-Abelian Seiberg Witten equations 
\begin{eqnarray*}
F_A^{+} &=&\tau \left( \phi \right)  \\
\frak{D}_A\phi  &=&0,
\end{eqnarray*}
on $E$ has a solution $\left( D_A,\phi \right) $, then $\Bbb{D}_{\left(
A,\phi \right) }=\left( 
\begin{array}{ll}
\nabla _L^{LC} & \tilde{\phi} \\ 
\widetilde{\phi }^{*} & D_A
\end{array}
\right) $ is a perturbed anti-self-dual $U\left( 2,2\right) $ connection on $%
\Bbb{E=}S_L^{-}\bigoplus E.$
\end{proposition}

Perturbed anti-self-dual bundle satisfies certain topological constraint,
namely the Chern number inequality. When the equality holds, the connection
under considerations is in fact a projectively flat connection.

\begin{proposition}
Let $\Bbb{E}$ be any complex vector bundle of rank $r$ over a Riemannian
four manifold $X.$ Suppose that $\Bbb{E}$ has a perturbed anti-self-dual $%
U\left( p,q\right) $ connection, then we have 
\[
c_1^2\left( E\right) \leq \frac{2r}{r-1}c_2\left( E\right) .
\]

Moreover equality sign holds if and only if $\Bbb{E}$ is a projectively flat 
$U\left( p,q\right) $ bundle.
\end{proposition}

The proof is the same as the Chern number inequality in the theory of
Hermitian-Einstein-Higgs bundle \cite{S} and this proposition is also stated
in \cite{L}.

Now if we take $\Bbb{E}$ to be $S_{L}^{-}\bigoplus S_{L}^{+}$, then the
above Chern number inequality is equivalent to 
\[
\sigma \left( X\right) \leq 0. 
\]
Notice that the Chern class of $L$ does not show up in this inequality.
Combining these two propositions we have the following theorem:

\begin{theorem}
Let $X$ be an Einstein four manifold with $Spin^c$ structure $L$ on $X$.
Suppose that the non-Abelian Seiberg Witten equations 
\begin{eqnarray*}
F_A^{+}=\tau \left( \phi \right)  \\
\frak{D}_A\phi =0,
\end{eqnarray*}
on $E=S_L^{+}$ has a solution $\left( D_A,\phi \right) $, then $\Bbb{D}%
_{\left( A,\phi \right) }=\left( 
\begin{array}{ll}
\nabla _L^{LC} & \tilde{\phi} \\ 
\widetilde{\phi }^{*} & D_{A,L}
\end{array}
\right) $ is a perturbed anti-self-dual $U\left( 2,2\right) $ connection on $%
\Bbb{E=}S_L^{-}\bigoplus E.$

In particular we have characteristic number inequality 
\[
\sigma \left( X\right) \leq 0,
\]
Moreover, equality sign holds if and only if $\Bbb{D}_{\left( A,\phi \right)
}$ is a projectively flat $U\left( 2,2\right) $ connection over $X.$
\end{theorem}

\section{Proof of the main theorem}

In this section we prove our main theorem:

\begin{theorem}
Let $X$ be an Einstein four manifold with $\sigma \left( X\right) =0.$
Suppose that for some $Spin^c$ structure $L$ on $X$, the non-Abelian Seiberg
Witten equations 
\begin{eqnarray*}
F_A^{+} &=&\tau \left( \phi \right)  \\
\frak{D}_A\phi  &=&0,
\end{eqnarray*}
on $E=S_L^{+}$ has a solution, then the universal covering of $X$ must be $%
\Bbb{R}^4,B^2\times B^2$ or $B^4$ with a locally symmetric Riemannian metric.
\end{theorem}

Proof of theorem: From the previous section we know that $\left( D_A,\phi
\right) $ induces a projectively flat connection on $S_L^{-}\bigoplus
E=S_L^{-}\bigoplus S_L^{+}.$

From the flatness equation we have $D_A\tilde{\phi}=0.$ It follows that $%
\phi $ is a parallel spinor, $D_A\phi =0.$ When regarding $\phi $ as a
homomorphism from $S_L^{+}$ to $E,$ the rank of $Ker\left( \phi \right) $ is
constant and $Ker\left( \phi \right) $ is a vector bundle over $M$. We
divide into the following three cases : (i) $rank\left( Ker\phi \right) =2$
(ii) $rank\left( Ker\phi \right) =0$ and (iii) $rank\left( Ker\phi \right)
=1.$

(i) First case: $rank\left( Ker\phi \right) =2,$ that is $\phi $ is a
trivial homomorphism. Hence the second fundamental form of the connection $%
\Bbb{D}_{\left( A,\phi \right) }$ on $S_L^{-}\bigoplus E$ is zero. Now the
projectively flatness of $S_L^{-}\bigoplus E$ implies that $S_L^{-}$ is a
projectively flat $U\left( 2\right) $ bundle with respect to the Levi-Civita
connection twisted by $D_L$.

Using curvature decomposition for $X$ and its associated bundles, we have $%
W^{-}=0=s$. Applying the equality 
\[
\sigma \left( X\right) =\frac 1{12\pi ^2}\int_X\left( \left| W^{+}\right|
^2-\left| W^{-}\right| ^2\right) dv, 
\]
we obtain $W^{+}=0$ because $\sigma \left( X\right) =0.$ Together with the
fact that $Rc^0=0$ we know that $X$ has vanishing Riemann curvature tensor.
That is the universal covering of $X$ is $\Bbb{R}^4$ with the flat metric,
in particular $\chi \left( X\right) =0.$

(ii) Second case: $Ker\left( \phi \right) =\left\{ 0\right\} .$ That is $%
\phi $ is a parallel isomorphism 
\[
\phi :S_L^{+}\stackrel{\cong }{\rightarrow }E_L. 
\]

Consider $g=\left( \phi \phi ^{*}\right) ^{-1/2}\phi $. Then $g$ is a
unitary gauge transformation of $E.$ Moreover it carries $\nabla _{L}^{LC}$
to $D_{A},$ namely 
\[
D_{A}=g\circ \nabla _{L}^{LC}\circ g^{-1}. 
\]

Therefore by composing $E$ with a unitary transformation if necessary, we
can assume that $D_A=\nabla _L^{LC}$ and $g$ is the identity transformation.
From the projective flatness for the twisted Levi-Civita connection on $%
S_L^{-}\bigoplus S_L^{+}$ and the curvature decomposition for four
manifolds, we know that $W^{+},W^{-}$ and $s$ are all parallel. Together
with the vanishing of $Rc^0,$ the whole Riemann curvature tensor is
parallel: $\nabla \left( Rm\right) =0.$ This means that $g$ is a locally
symmetric metric. We now determine the locally symmetric space explicitly.

Since $D_A=\nabla _L^{LC},$ the trace free part of $\phi $ is a parallel
self-dual (complex) two form on $X.$ If this self-dual two form is zero,
then $\phi =\frac 12\left( Tr\phi \right) I_{S^{+}}$ where $Tr\phi $ is a
constant function on $X.$ Using the identification $\Lambda _{-}^2=\mathbf{su%
}\left( S_L^{-}\right) $ we have 
\[
\Omega _{-}^2\left( X,\mathbf{u}\left( S_L^{-}\right) \right) =\Omega
_{-}^2\left( X,\mathbf{su}\left( S_L^{-}\right) \bigoplus \Bbb{R}\right)
=\Gamma \left( X,End\left( \Lambda _{-}^2\right) \bigoplus \Lambda
_{-}^2\right) . 
\]
It is easy to check that the component of $P_{-}\left( \tilde{\phi}^{*}%
\tilde{\phi}\right) $ in $\Gamma \left( X,End\left( \Lambda _{-}^2\right)
\right) \subset $ $\Gamma \left( X,End\left( \Lambda _{-}^2\right) \bigoplus
\Lambda _{-}^2\right) =\Omega _{-}^2\left( X,\mathbf{u}\left( S_L^{-}\right)
\right) $ equals $cI_{\Lambda _{-}^2}$ for some constant $c$ (morover the
component of $P_{-}\left( \tilde{\phi}^{*}\tilde{\phi}\right) $ in $\Gamma
\left( X,\Lambda _{-}^2\right) $ is proportional to $F_L^{-},$ even though
we will not use this fact).

Applying the flatness condition, we have the trace free part of $%
F_{S_L^{-}}^{-}$ equals $cI_{\Lambda _{-}^2}$. This implies $W^{-}=0$ and
the scalar curvature $s$ of $X$ is a constant function. Using 
\[
0=\sigma \left( X\right) =\frac 1{12\pi ^2}\int_X\left( \left| W^{+}\right|
^2-\left| W^{-}\right| ^2\right) dv, 
\]
we also have $W^{+}=0.$ As a result $X$ has constant sectional curvature. By
checking the sign of $c$ we know $s$ is negative and $X$ has constant
negative sectional curvature. Hence $X=\Gamma \backslash SO^0\left(
4,1\right) /SO\left( 4\right) =\Gamma \backslash B^4.$

Now if trace free part of $\phi $ is nonzero, then we have a nontrivial
parallel self-dual two form on $X.$ As a result, the holonomy group of $%
\left( X,g\right) $ will be reduced to $U\left( 2\right) $. That is $g$ is a
K\"{a}hler metric on $X$ and the above parallel self-dual two form is
proportional to its K\"{a}hler form $\omega $. In particular $\left(
X,g\right) $ is a K\"{a}hler Einstein surfaces.

Recall that the Levi-Civita connection on $X$ twisted by $D_L$ induces a
projectively flat connection on $S_L^{-}\bigoplus S_L^{+}.$ Reversal of
orientation will interchange $S_L^{+}$ and $S_L^{-},$ but their direct sum $%
S_L^{-}\bigoplus S_L^{+}$ remains unchanged. Moreover we just showed that
the projectively flat connection on $S_L^{-}\bigoplus S_L^{+}$ is induced
from the Levi-Civita connection twisted by $D_L$ and therefore also remain
unchanged after the reversal of orientations. Unlike the anti-self-duality
condition, flatness is insensitive to the orientation of the manifold.

Changing the orientation of $X,$ the Riemannian metric on $X$ remain
Einstein because the metric tensor is independent of the orientation on $X$.
Even though the signature will change sign under the reversal of
orientation, we still have $\sigma \left( X\right) =0.$ Repeat the same
arguments as above, we show that $X$ is K\"{a}hler Einstein with respect to
both orientations. Therefore the Seiberg Witten invariant for $X$ with both
orientations are non-trivial. By the result of \cite{L}, we concludes that $%
X=\Gamma \backslash B^2\times B^2.$

(iii) Third case: $rank\left( Ker\phi \right) =1.$ We shall see that this
case cannot happen. Let $L_1=Ker\phi \subset S_L^{+}$ and $L_2=\left(
L_1\right) ^{\bot }\subset S_L^{+}.$ Then we have $L_2=\mbox{Im}\phi \ $and
an orthogonal decomposition $S_L^{+}=L_1\bigoplus L_2.$ Since $\phi $ is
parallel and $E$ and $S_L^{+}$ are isometric to each other, we also have an
orthogonal decomposition of $E$: $E=L_1\bigoplus L_2.$ With respect these
decompositions, we have 
\[
\phi =\left( 
\begin{array}{ll}
0 & 0 \\ 
0 & \beta
\end{array}
\right) , 
\]
for some non-zero constant $\beta $. Therefore $Tr_E\tau \left( \phi \right)
\in \Omega _{+}^2$ is a \textit{nonzero} parallel self-dual two form on $X.$
This implies that the holonomy group of $\left( X,g\right) $ is reduced to $%
U\left( 2\right) $ or equivalently $g$ is a K\"{a}hler metric on $X$ with
respect to some integrable complex structure. If we denote the K\"{a}hler
form by $\omega ,$ then it is proportional to $Tr_E\tau \left( \phi \right)
. $

Then the component of $P_{-}\left( \tilde{\phi}^{*}\tilde{\phi}\right) $ in $%
\Gamma \left( X,End_0\left( \Lambda _{-}^2\right) \right) \subset \Omega
_{-}^2\left( X,\mathbf{u}\left( S_L^{-}\right) \right) $ is zero by direct
computation. Using the projectively flatness condition and curvature
decomposition, we have $W^{-}=0$ and in fact the universal covering of $X$
must be $B_{\Bbb{C}}^2$ \cite{L}. In particular we have the Chern number
equality $c_1^2\left( X\right) =3c_2\left( X\right) >0.$ This violate our
assumption that $\sigma \left( X\right) =0$ because $\sigma \left( X\right)
=\left( c_1^2\left( X\right) -2c_2\left( X\right) \right) /3.$ Therefore $%
rank\left( Ker\phi \right) $ can never be one.$\square $

Combining these results and the work of Yau \cite{Y} and LeBrun \cite{LB},
we obtain the following characterization of four dimensional locally
symmetric spaces of noncompact type.

\begin{theorem}
Let $X$ be an Einstein four manifold with $\chi \left( X\right) >0$. Then $X$
is a locally symmetric space of noncompact type iff one of the followings
holds

(1) the non-Abelian Seiberg Witten equation for $S_L^{+}$ has a solution and 
$\sigma \left( X\right) =0$.

(2) the Seiberg Witten equation for $L$ has a solution and $3\sigma \left(
X\right) =\chi \left( X\right) .$

The universal cover of $X$ is isometric to $B^2\times B^2$ or $B^4$ in (1)
and $B_{\Bbb{C}}^2$ in (2).
\end{theorem}

Notice that when $\tilde{X}=B_{\Bbb{C}}^2,$ Seiberg Witten equation has a
unique solution provided by the K\"{a}hler Einstein metric on $X.$ Moreover
the number of solutions to the Seiberg Witten equation counted with
multiplicity is a differentiable invariant. As a result, as LeBrun pointed
out in \cite{LB}, Einstein metric on such manifolds is unique. Therefore an
important question is to use the space of solutions to the non-Abelian
Seiberg Witten equations to define Donaldson type differentiable invariant
for $X$. The author is informed by P.Feehan that this might be possible in
some situations, modulo certain technical difficulties. A positive solution
to the above question will tell us that the above characterization of
locally symmetric space depends only on Einstein metric and the
differentiable structure of $X.$ In particular, Einstein metric on locally
symmetric four manifolds of noncompact type would be unique.
\begin{verbatim}
Acknowledgement: The author would like to thank S.Adams, J.Bryan, 
P.Feehan,D.Kotschick, A.K.Liu, P.Lu and Y.S.Poon for discussions. 
The author is also grateful to S.T.Yau who taught me the theory 
of uniformizations.
 
\end{verbatim}

\end{document}